\documentclass[11pt]{article}

\usepackage{amssymb}
\usepackage{latexsym}

\usepackage{graphicx}

\newtheorem{defi}{Definition}[section]
\newtheorem{theo}{Theorem}[section]

\newtheorem{remark}[theo]{Remark}

\newtheorem{lemma}[theo]{Lemma}

\newtheorem{con}[theo]{Conjecture}

\newtheorem{fact}[theo]{Fact}

\newcommand{\qed}{\hspace*{\fill} \rule{7pt}{7pt}}

\begin{document}

\title{On Lagrangians of $r$-uniform Hypergraphs}

\author{ Yuejian Peng \thanks{ Supported by National Natural Science Foundation of China (No. 11271116). School of Mathematics, Hunan University, Changsha, 410082, P.R. China. Email: ypeng1@163.com} \and Qingsong Tang \thanks{Mathematics School, Institute of Jilin University, Changchun, 130012, P.R. China, and College of Sciences, Northeastern University, Shenyang, 110819, P.R. China. 
t\_qsong@sina.com.cn} \and  Cheng Zhao \thanks{School of Mathematics, Jilin University, Changchun 130022, P.R. China and Department of Mathematics and Computer Science, Indiana State University, Terre Haute, IN, 47809. Email: cheng.zhao@indstate.edu}
}


\maketitle
\date

\begin{abstract}
 A remarkable connection between the order of a maximum clique and the Lagrangian of a graph was established by Motzkin and Straus in \cite{MS}.  This connection and its extensions were successfully employed in optimization to provide heuristics for the maximum clique number in graphs. It  has been also applied in spectral graph theory.  Estimating the Lagrangians of hypergraphs has been  successfully applied in the course of studying the Tur\'an densities of several hypergraphs  as well. It is useful in practice if Motzkin-Straus type results hold for hypergraphs. However, the obvious generalization of Motzkin and Straus' result to hypergraphs is false. We attempt to explore the relationship between the Lagrangian of a hypergraph and the order of its maximum cliques for hypergraphs when the number of edges is in certain range. In this paper, we give some  Motzkin-Straus type results for $r$-uniform hypergraphs. These results generalize  and refine a result of Talbot in \cite{T} and a result in \cite{PZ}. 
\end{abstract}

keywords: Cliques of Hypergraphs; Lagrangians of Hypergraphs; Optimization.

\section{Introduction}
In 1965, Motzkin and Straus \cite{MS} established  a continuous characterization of the clique number of a graph using the Lagrangian of a graph. Namely, the Lagrangian of a graph is the Lagrangian of its maximum clique which is determined by the order of a maximum cliques.  Applying this connection, they provided a new proof of  classical Tur\'an's theorem \cite{Turan} on the extremal number of a complete graph. This connection has  been also applied in spectral graph theory \cite{wilf}.
 Furthermore, the Motzkin-Straus result and its extension were successfully employed in optimization to provide heuristics for the maximum clique problem.  The Motzkin-Straus theorem has  been also generalized to vertex-weighted graphs \cite{G9} and edge-weighted graphs with applications to pattern recognition in image analysis (see 
\cite{B2}, \cite{B3}, \cite{G9}, \cite{PP}, \cite{PP15}, \cite{RTP20}).   It is interesting to explore whether  similar results hold for hypergraphs. 
The obvious generalization of Motzkin and Straus' result to hypergraphs is false. In fact, there are many examples of hypergraphs that do not achieve their Lagrangian on any proper subhypergraph.
In this paper, we provide evidences that the Lagrangian of an $r$-uniform hypergraph is related to the order of its maximum cliques under some conditions. Some definitions and notations are needed in order to state the questions and results precisely.

  Let ${\mathbb N}$ be the set of all positive integers. Let  $V$ be a set and $r\in {\mathbb N}$.  Let  $V^{(r)}$ denote the family of all $r$-subsets of $V$. An {\em $r$-uniform graph} or {\em $r$-graph} $G$ is a set $V(G)$ of vertices together with a set $E(G) \subseteq V(G) ^{(r)}$ of edges. An edge $e=\{a_1, a_2, \ldots, a_r\}$ will be simply denoted by $a_1a_2 \ldots a_r$. An $r$-graph $H$ is  a {\it subgraph} of an $r$-graph $G$, denoted by $H\subseteq G$ if $V(H)\subseteq V(G)$ and $E(H)\subseteq E(G)$.    Let $K^{(r)}_t$ denote the complete $r$-graph on $t$ vertices, that is the $r$-graph on $t$ vertices containing all possible edges. A complete $r$-graph on $t$ vertices is also called a clique with order $t$. For  $n \in {\mathbb N}$, we denote the set $\{1, 2, 3, \ldots, n\}$ by $[n]$. Let $[n]^{(r)}$  represent the  complete $r$-uniform graph on the vertex set $[n]$. When $r=2$, an $r$-uniform graph is a simple graph.  When $r\ge 3$,  an $r$-graph is often called a hypergraph.

\begin{defi}
Let $G$ be  an $r$-uniform graph  with vertex set $[n]$ and
edge set $E(G)$.  Let $S=\{\vec{x}=(x_1,x_2,\ldots ,x_n)\in R^n: \sum_{i=1}^{n} x_i =1, x_i
\ge 0 {\rm \ for \ } i=1,2,\ldots , n \}$. For $\vec{x}=(x_1,x_2,\ldots ,x_n)\in S$,
define
$$\lambda (G,\vec{x})=\sum_{i_1i_2 \cdots i_r \in E(G)}x_{i_1}x_{i_2}\ldots x_{i_r}.$$

 The Lagrangian of
$G$, denoted by $\lambda (G)$, is defined as $$\lambda (G) = \max \{\lambda (G, \vec{x}): \vec{x} \in S \}.$$

 A vector $\vec{y}\in S$ is called an optimal weighting for $G$ if $\lambda (G, \vec{y})=\lambda(G)$.
\end{defi}

 The following fact is easily implied by the definition of the Lagrangian.

\begin{fact}\label{mono}
Let $G_1$, $G_2$ be $r$-uniform graphs and $G_1\subseteq G_2$. Then $\lambda (G_1) \le \lambda (G_2).$
\end{fact}

 The following theorem by Motzkin and Straus in \cite{MS} shows that  the Lagrangian of a 2-graph is determined by the order of its maximum clique.

\begin{theo} (Motzkin and Straus \cite{MS}) \label{MStheo}
If $G$ is a 2-graph in which a largest clique has order $l$, then
$\lambda(G)=\lambda(K^{(2)}_l)=\lambda([l]^{(2)})={1 \over 2}(1 - {1 \over l})$.
\end{theo}

As mentioned earlier, there are many examples of hypergraphs that do not achieve their Lagrangian on any proper subhypergraph and the obvious generalization of Motzkin and Straus' result to hypergraphs is false.
S\'os and Straus attempted to generalize the Motzkin-Straus theorem to hypergraphs in \cite{SS}.  Recently,  Rota Bul\'o and Pelillo generalized the Motzkin and Straus' result to $r$-graphs in some way using a continuous characterization of maximal cliques other than Lagrangians of hypergraphs in \cite{BP1} and \cite{BP2}. Lagrangians of hypergraphs has been proved to be a useful tool in hypergraph extremal problems. For example, Frankl and R\"odl \cite{FR84} applied it in disproving Erd\"os' long standing jumping constant conjecture. It has also been applied in finding Tur\'an densities of hypergraphs in \cite{FF}, \cite{sidorenko89} and \cite{mubayi06}. We attempt to explore the relationship between the Lagrangian of a hypergraph and the  order of its maximum cliques for hypergraphs when the number of edges is in certain range though the obvious generalization of Motzkin and Straus' result to hypergraphs is false.  The following two conjectures are proposed in \cite{PZ}.

\begin{con} \label{conjecture1} (Peng-Zhao \cite{PZ})
Let $l$ and $m$ be positive integers satisfying ${l-1 \choose r} \le m \le {l-1 \choose r} + {l-2 \choose r-1}$.
Let $G$ be an $r$-graph with $m$ edges and  contain a clique of order  $l-1$. Then $\lambda(G)=\lambda([l-1]^{(r)})$.
\end{con}

The upper bound ${l-1 \choose r} + {l-2 \choose r-1}$ in this conjecture is the best possible. For example, if $m ={l-1 \choose r}+{l-2 \choose r-1}+1$ then $\lambda(C_{r,m}) > \lambda([l-1]^{(r)})$, where  $C_{r,m}$ is the $r$-graph on the  vertex set $[l]$ and with the edge set $[l-1]^{(r)}\cup\{i_1\cdots i_{r-1}l, i_1\cdots i_{r-1}\in [l-2]^{(r-1)}\}\cup\{1\cdots (r-2)(l-1)l\}$. Take  $\vec{x}=(x_1,\ldots,x_l) \in S$, where $x_1=x_2=\cdots=x_{l-2}={1 \over l-1}$ and $x_{l-1}=x_{l}={1 \over 2(l-1)}$.
Then $\lambda(C_{r,m})\ge \lambda(C_{r,m},\vec{x} )> \lambda([l-1]^{(r)})$.

\begin{con} \label{conjecture2} (Peng-Zhao \cite{PZ})
Let $l$ and $m$ be positive integers satisfying ${l-1 \choose r} \le m \le {l-1 \choose r} + {l-2 \choose r-1}$.
 Let $G$ be an $r$-graph with $m$ edges and contain no clique of order $l-1$.  Then $\lambda(G) < \lambda([l-1]^{(r)})$.
\end{con}

In \cite{PZ}, Conjecture \ref{conjecture1} is proved for $r=3$.

\begin{theo} \label{theorem 1} (Peng-Zhao \cite{PZ}) Let $m$ and $l$ be positive integers satisfying ${l-1 \choose 3} \le m \le {l-1 \choose 3} + {l-2 \choose 2}$. Let $G$ be a $3$-graph with $m$ edges and $G$ contain a clique of order  $l-1$. Then $\lambda(G) = \lambda([l-1]^{(3)})$.
\end{theo}

In \cite{PZZ}, an algorithm is proposed to check the validity of Conjecture \ref{conjecture2} for 3-graphs and, as a demonstration, that algorithm confirms Conjecture \ref{conjecture2} for some small $l$. For 3-graphs, the validity of Conjecture \ref{conjecture2} for some small $l$ is verified in \cite{PZ2} as well.


In \cite{FF}, Frankl and F\"uredi applied the Lagrangians of related hypergraphs to estimate Tur\'an densities of hypergraphs. They asked  the following question: Given $r \ge 3$ and $m \in {\mathbb N}$ how large can the Lagrangian of an $r$-graph with $m$ edges be? An answer to the above question would be quite useful in estimating  Tur\'an densities of  hypergraphs.

The following definition is needed in order to state their conjecture on this problem. For distinct $A, B \in {\mathbb N}^{(r)}$,  $A$ is less than $B$ in the {\em colex ordering} if $max(A \triangle B) \in B$, where $A \triangle B=(A \setminus B)\cup (B \setminus A)$. For example, $246 < 156$ in ${\mathbb N}^{(3)}$ since $max(\{2,4,6\} \triangle \{1,5,6\}) \in \{1,5,6\}$. In colex ordering, $123<124<134<234<125<135<235<145<245<345<126<136<236<146<246<346<156<256<356<456<127<\cdots .$ Note that the first $l \choose r$ $r$-tuples in the colex ordering of ${\mathbb N}^{(r)}$ are the edges of $[l]^{(r)}$. Let $C_{r,m}$ denote the $r$-graph with $m$ edges formed by taking the first $m$ elements in the colex ordering of ${\mathbb N}^{(r)}$. 

\begin{con} (Frankl and F\"uredi \cite{FF})\label{conjecture} The $r$-graph with $m$ edges formed by taking the first $m$ sets in the colex ordering of ${\mathbb N}^{(r)}$ has the largest Lagrangian of all $r$-graphs with $m$ edges. In particular, the $r$-graph with $l \choose r$ edges and the largest Lagrangian is $[l]^{(r)}$.
\end{con}

 Theorem \ref{MStheo} implies that this conjecture is true when $r=2$. For the case $r=3$, Talbot in \cite{T} proved the following result.

\begin{theo} (Talbot \cite{T}) \label{Tal} Let $m$ and $l$ be positive integers satisfying
$${l-1 \choose 3} \le m \le {l-1 \choose 3} + {l-2 \choose 2} - (l-1).$$
Then Conjecture \ref{conjecture} is true for $r=3$ and this value of $m$.  Conjecture \ref{conjecture} is also true for $r=3$ and $m= {l \choose 3}-1$ or $m={l \choose 3} -2$.
\end{theo}

For the case $r=3$, Tang, Peng, Zhang, and Zhao in \cite{TPZZ1} and \cite{TPZZ2} proved the following.
\begin{theo} (Tang, Peng, Zhang, and Zhao \cite{TPZZ1},\cite{TPZZ2} ) \label{Tang} Let $m$ and $t$ be integers.
Then  Conjecture \ref{conjecture} is  true for $r=3$ and $m= {t \choose 3}-3$, $m={t \choose 3} -4$   or $m={t \choose 3} -5$.
\end{theo}

In \cite{HPZ}, He, Peng, and Zhao verified Frankl and F\"uredi's conjecture  for $m\le 50$ when $r=3$.

The truth of Frankl and F\"uredi's conjecture is not known in general for $r \ge 4$. Even in the case $r=3$, it is still open.

The following  result was given in \cite{T}.

\begin{lemma} \cite{T} \label{LemmaTal7}
For positive integers $m,l,$ and $r$ satisfying ${l-1 \choose r} \le m \le {l-1 \choose r} + {l-2 \choose r-1}$,
we have $\lambda(C_{r,m}) = \lambda([l-1]^{(r)})$.
\end{lemma}

If Conjectures \ref{conjecture1} and \ref{conjecture2} are true, then Conjecture \ref{conjecture} is true for this range of $m$. In Section \ref{sectionconjecture1}, we provide some evidences for  Conjectures \ref{conjecture1} and \ref{conjecture2}. 
In addition to several other results,  we  will prove the following result in Section \ref{sectionconjecture1}.

\begin{theo} \label{Corollary 2}
 (a) Let $m$ and $l$ be positive integers satisfying ${l-1 \choose r} \le m \le {l-1 \choose r} + {l-2 \choose r-1}-(2^{r-3}-1)({l-2 \choose r-2}-1)$.  Let $G$ be an $r$-graph on $l$ vertices with $m$ edges and contain a clique of order $l-1$. Then $\lambda(G) = \lambda([l-1]^{(r)})$. 


(b) Let $m$ and $l$ be positive integers satisfying ${l-1 \choose 3} \le m \le {l-1 \choose 3} + {l-2 \choose 2}-(l-2)$.  Let $G$ be a $3$-graph with $m$ edges and without containing a clique of order $l-1$. Then $\lambda(G) <\lambda([l-1]^{(3)})$. 
\end{theo} 

When $r=3$, Theorem \ref{Corollary 2} (a) and Lemma \ref{lemmaclique1}  imply Theorem \ref{theorem 1}.   Theorem \ref{Corollary 2} (b) and Theorem \ref{theorem 1} refine Theorem \ref{Tal}.  

Theorem \ref{Corollary 2} gives a connection between a continuous optimization problem  and the maximum clique problem of 3-uniform hypergraphs. Since  practical problems such as computer vision and image analysis are related to the maximum clique problems in hypergraphs, this type of results open a  door to  such practical applications.  The  results in this paper can be applied in estimating Lagrangians of some hypergraphs, for example, calculations involving estimating Lagrangians of several hypergraphs in \cite{FF} can be much simplified when applying the  results in this paper.

Some preliminary results  will be stated in the following section.

\section{Preliminary Results}
For an $r$-graph $G=(V,E)$ on the vertex set $[n]$ and $i \in V$, let $E_i=\{A \in V^{(r-1)}: A \cup \{i\} \in E\}$ be the link of the vertex $i$. Similarly, for a pair of vertices $i,j \in V$, let $E_{ij}=\{B \in V^{(r-2)}: B \cup \{i,j\} \in E\}$. Let $E^c_i=\{A \in V^{(r-1)}: A \cup \{i\} \in V^{(r)} \backslash E\}$, and 
$E^c_{ij}=\{B \in V^{(r-2)}: B \cup \{i,j\} \in V^{(r)} \backslash E\}$. Denote $$E_{i\setminus j}=E_i\cap E^c_j.$$  Let us  impose one additional condition on any optimal weighting ${\vec x}=(x_1, x_2, \ldots, x_n)$ for an $r$-graph $G$:
\begin{eqnarray}
  &&|\{i : x_i > 0 \}|{\rm \ is \ minimal, i.e. \ if}  \ \vec y \in S
  {\rm \ satisfies } \ \  |\{i : y_i > 0 \}| < |\{i : x_i > 0 \}|, \nonumber \\
&&{\rm \  then \ } \lambda (G, {\vec y}) < \lambda(G) \label{conditionb}.
\end{eqnarray}

Note that $\lambda (E_i, {\vec x})$ corresponds to the partial derivative of  $\lambda(G, \vec x)$ with respect to $x_i$.
The following lemma gives some necessary conditions of an optimal weighting of  $\lambda(G)$.

\begin{lemma} (Frankl and R\"odl \cite{FR84}) \label{LemmaTal5} Let $G=(V,E)$ be an $r$-graph on the vertex set $[n]$ and ${\vec x}=(x_1, x_2, \ldots, x_n)$ be an optimal weighting for $G$ with $k$  ($\le n$) non-zero weights satisfying condition (\ref{conditionb}). Then for every $\{i, j\} \in [k]^{(2)}$, (a) $\lambda (E_i, {\vec x})=\lambda (E_j, \vec{x})=r\lambda(G)$, (b) there is an edge in $E$ containing both $i$ and $j$.
\end{lemma}

\begin{defi}
An $r$-graph $G=(V,E)$ on the vertex set $[n]$ is left-compressed if $j_1j_2\ldots j_r\in E$ implies $i_1i_2\ldots i_r\in E$ whenever $i_k \leq j_k, 1\leq k \leq r$. Equivalently, an $r$-graph $G=(V,E)$  on the vertex set $[n]$ is  left compressed if $E_{j\setminus i}=\emptyset$ for any $1\le i<j\le n$.
\end{defi}


\begin{remark}\label{r1} (a) In Lemma \ref{LemmaTal5}, part (a) implies that
$$x_j\lambda(E_{ij}, {\vec x})+\lambda (E_{i\setminus j}, {\vec x})=x_i\lambda(E_{ij}, {\vec x})+\lambda (E_{j\setminus i}, {\vec x}).$$
In particular, if $G$ is left compressed, then
\begin{equation}\label{enbhd}
(x_i-x_j)\lambda(E_{ij}, {\vec x})=\lambda (E_{i\setminus j}, {\vec x})
\end{equation}
for any $i, j$ satisfying $1\le i<j\le k$ since $E_{j\setminus i}=\emptyset$.

(b) By (\ref{enbhd}), if  $G$ is left-compressed, then an optimal weighting  ${\vec x}=(x_1, x_2, \ldots, x_n)$ for $G$  must satisfy
\begin{equation}\label{conditiona}
x_1 \ge x_2 \ge \ldots \ge x_n \ge 0.
\end{equation}
\end{remark}

Denote 
\begin{eqnarray*}
&&\lambda_m^r=\max\{\lambda(G): G {\rm \ is \ an \ } r-{\rm graph \ with \ } m {\rm \ edges }\}, \\
&&\lambda_{(m,l)}^r=\max\{\lambda(G): G {\rm \ is \ an \ } r-{\rm graph \ with \ } m {\rm \ edges \ and  \ contains \ a \ clique \ of \ order\ } l  \}, 
{\rm \ and } \\
&&\lambda_{(m,l)}^{r-}=\max\{\lambda(G): G {\rm \ is \ an \ } r-{\rm graph \ with \ } m {\rm \ edges \ and  \ without \  a \ clique \ of \ order\ } l \}. 
\end{eqnarray*}

The following two lemmas imply that we only need to consider left-compressed $r$-graphs when Conjecture \ref{conjecture} and Conjecture \ref{conjecture1} are explored.

\begin{lemma} \label{lemmaleftcompress}\cite{T} There exists a left compressed $r$-graph $G$ with $m$ edges  such that $\lambda(G)=\lambda_m^r$. 
\end{lemma}

\begin{lemma} \label{lemma1}\cite{PZ} Let $m$ and $l$  be  positive integers satisfying $m\le {l \choose r} -1$. Then there exists a left compressed $r$-graph $G$ containing the clique  $[l-1]^{(r)}$ with $m$ edges  such that $\lambda(G)=\lambda_{(m,l-1)}^r$.
\end{lemma}

When Conjectures \ref{conjecture1}  and \ref{conjecture2} were  discussed for $r=3$ in \cite{PZ} and \cite{PZZ}, the following results were proved.

\begin{lemma}\label{lemmaclique1}\cite{PZ} Let $m$ and $l$  be  positive integers satisfying ${l-1 \choose 3}\le m\le {l-1 \choose 3}+{l-2 \choose 2}$. Then there exists a left compressed $3$-graph $G$  on the vertex set  $[l]$ with $m$ edges and containing the clique  $[l-1]^{(3)}$  such that $\lambda(G)=\lambda_{(m,l-1)}^{3}$.
\end{lemma}

\begin{lemma}\label{lemmanoclique}\cite{PZZ} Let $m$ and $l$  be  positive integers satisfying ${l-1 \choose 3}\le m\le {l-1 \choose 3}+{l-2 \choose 2}$. Then there exists a left compressed $3$-graph $G$ on the vertex set  $[l]$ with $m$ edges and without containing the clique  $[l-1]^{(3)}$ such that $\lambda(G)=\lambda_{(m,l-1)}^{3-}$.
\end{lemma}

\section{Evidence for Conjectures  \ref{conjecture1} and \ref{conjecture2}}\label{sectionconjecture1}

Frank and F\"uredi \cite{FF} originally asked how large the Lagrangian of an $r$-graph with   $l$ vertices and $m$ edges can be, where $m \le {l \choose r}$. For a given $r$-graph with $l$ vertices and $m$ edges, let
$$\lambda(l,r,m)=max \{\lambda(G): \,\ G=(V,E)\,\ {\rm is \,\ an \,\ r-graph}, |V|=l,|E|=m\}.$$

In \cite{T}, the following result is proved, which is the evidence for Conjecture \ref{conjecture1} for $r$-graphs $G$ on exactly $l$ vertices.

\begin{theo} (Talbot \cite{T}) \label{Talr} For any $r \ge 4$ there exists constants $\gamma_r$ and $\kappa_0(r)$ such that if $m$ satisfies
$${l-1 \choose r} \le m \le {l-1 \choose r} + {l-2 \choose r-1} - \gamma_r (l-1)^{r-2},$$
with $l \ge \kappa_0(r)$, then $\lambda(l,r,m)=\lambda(C_{r,m})=\lambda([l-1]^{(r)})$.
\end{theo}

In \cite{TPZZ1}, we proved: 

\begin{theo} (Tang, Peng Zhang, and Zhao \cite{TPZZ1}) \label{theorem 2} Let $m$ and $l$ be positive integers  satisfying ${l \choose r}-4 \le m \le {l \choose r}-1$. Then the $r$-graph with $m$ edges formed by taking the first $m$ sets in the colex ordering of ${\mathbb N}^{(r)}$ has the largest Lagrangian of all $r$-graphs  with $m$ edges and 
$l$ vertices.
\end{theo}


Next, we point out a useful lemma.

\begin{lemma}\label{lemmaadd} Let $G$ be a left-compressed $r$-graph on the vertex set $[l]$ containing the clique $[l-1]^{(r)}$. Let ${\vec x}=(x_1, x_2, \ldots, x_l)$ be an optimal weighting for $G$. Then
\begin{equation}\label{eq234}
 x_1\le x_{l-1}+x_l \le 2x_{l-1}.
\end{equation}
\end{lemma}

{\em Proof.} Note that $x_{l-1}>0$. Otherwise $\lambda(G, {\vec x})\le \lambda([l-2]^{(r)})$ contradicting to that ${\vec x}$ is an optimal weighting for $G$. Since $G$ is left compressed, applying  Remark \ref{r1}(a) by taking $i=1$, $j=l-1$, we get
\begin{equation}\label{addeq1}
x_1=x_{l-1}+{\lambda(E_{1 \setminus (l-1)}, \vec{x}) \over \lambda(E_{1(l-1)}, \vec{x})}.
\end{equation}
Since $G$ contains the clique $[l-1]^{(r)}$, then any $(r-1)$-tuple in $E_{1 \setminus (l-1)}$ must contain $l$ but not $1$ or $l-1$. Therefore
\begin{equation}\label{addeq2}
\lambda(E_{1 \setminus (l-1)}, \vec{x})\le\sum_{2\le i_1<i_2<\cdots <i_{r-2}\le l-2} x_{i_1}x_{i_2}\cdots x_{i_{r-2}}x_l.
\end{equation}
Since $G$ contains the clique $[l-1]^{(r)}$, then every $(r-2)$-tuple in $\{2, 3, \cdots, l-2\}^{(r-2)}$ belongs to $E_{1(l-1)}$, then 
\begin{equation}\label{addeq3}
\lambda(E_{1(l-1)}, \vec{x})\ge \sum_{2\le i_1<i_2<\cdots <i_{r-2}\le l-2} x_{i_1}x_{i_2}\cdots x_{i_{r-2}}.
\end{equation}
Combining inequalities (\ref{addeq2}) and (\ref{addeq3}), we get
\begin{equation}\label{addeq4}
{\lambda(E_{1 \setminus (l-1)}, \vec{x}) \over \lambda(E_{1(l-1)}, \vec{x})} \le x_l.
\end{equation}
 Applying inequality (\ref{addeq4}) to (\ref{addeq1}), we get (\ref{eq234}). \qed

Note that the only left-compressed $r$-graph on the vertex set $[r+1]$ is $C_{r, r+1}$. So we assume that an $r$-graph has at least $r+2$ vertices in this paper.

Next we give some results  refining Theorem \ref{Tal} when $r=3$. 

\begin{theo} \label{theorem 2a} Let $r\ge 3$ and $l\ge r+2$ be positive integers.
Let $G$ be a left-compressed $r$-graph  on the vertex set $[l]$  satisfying $|[l-2]^{(r-1)} \backslash E_l| \ge 2^{r-3}|E_{(l-1)l}|$.

(a) If $G$ contains $[l-1]^{(r)}$, then  $\lambda(G) = \lambda([l-1]^{(r)})$.

(b) If $G$  does not contain $[l-1]^{(r)}$, then  $\lambda(G) <\lambda([l-1]^{(r)})$.

 \end{theo}

{\em Proof.}  
(a) If $G$ contains $[l-1]^{(r)}$, then  clearly $\lambda(G)\ge  \lambda([l-1]^{(r)})$. We show that $\lambda(G) \le \lambda([l-1]^{(r)})$ as well.
Let  ${\vec x}=(x_1, x_2, \ldots, x_l)$ be an optimal weighting for $G$. Since $G$ is left-compressed, by Remark \ref{r1}(a), $x_1\ge x_2 \ge \cdots \ge x_l \ge 0$. If $x_{l}=0$, then the conclusion holds obviously, so we assume that $x_{l}>0$.

Consider a new weighting for $G$, ${\vec z}=(z_1, z_2, \ldots, z_l)$ given by $z_i=x_i$ for $i\neq l-1, l$, $z_{l-1}=0$ and $z_l=x_{l-1}+x_l$. By Lemma \ref{LemmaTal5}(a), $\lambda(E_{l-1}, \vec{x})=\lambda(E_{l}, \vec{x})$, so
\begin{eqnarray}\label{eq10b}
\lambda(G,\vec {z})- \lambda(G,\vec {x})&=&x_{l-1}(\lambda(E_{l}, \vec{x})-x_{l-1}\lambda(E_{(l-1)l}, \vec{x})) \nonumber \\
&&-x_{l-1}(\lambda(E_{l-1}, \vec{x})-x_{l}\lambda(E_{(l-1)l}, \vec{x}))-x_{l-1}x_l\lambda(E_{(l-1)l}, \vec{x}))\nonumber \\
&=&x_{l-1}(\lambda(E_{l}, \vec{x})-\lambda(E_{l-1}, \vec{x}))-x_{l-1}^2\lambda(E_{(l-1)l}, \vec{x})  \nonumber\\
&=& -x_{l-1}^2\lambda(E_{(l-1)l}, \vec{x}).
\end{eqnarray}


We will show that there exists a set of edges $F\subset \{1, ..., l-2, l\}^{(r)}\setminus E$ satisfying 
\begin{equation}\label{eq11b}
\lambda(F,\vec {z})\ge x_{l-1}^2\lambda(E_{(l-1)l}, \vec{x}).
\end{equation}

Then using (\ref{eq10b}) and (\ref{eq11b}), the $r$-graph $G^{*}=([l], E^{*})$, where $E^{*}=E\cup F$, satisfies $\lambda(G^{*}, \vec {z}))\ge \lambda(G)$. Since $\vec {z}$ has only $l-1$ positive weights, then $\lambda(G^{*}, \vec {z}))\le \lambda([l-1]^{(r)})$, and consequently,
$\lambda(G)\le \lambda([l-1]^{(r)}).$

We now construct the set of edges $F$. Let $D=[l-2]^{(r-1)} \setminus E_l$. Then by the assumption, 
$\vert D\vert \ge 2^{r-3}|E_{(l-1)l}|$ and 
$\lambda(D, \vec{x})\ge 2^{r-3}|E_{(l-1)l}|(x_{l-1})^{r-1}.$

Let $F$ consist of those edges in $\{1, ..., l-2, l\}^{(r)}\setminus E$ containing the vertex $l$. Then 
\begin{eqnarray}\label{eq132ab}
\lambda(F,\vec {z}) &=&(x_{l-1}+x_l)\lambda(D, \vec{x}) \nonumber  \\ 
&\ge& x_1 2^{r-3}|E_{(l-1)l}|(x_{l-1})^{r-1} \nonumber \\
&=& x_{l-1}^2 |E_{(l-1)l}|x_1( 2^{r-3}(x_{l-1})^{r-3}) \nonumber \\ 
&\ge&  x_{l-1}^2 |E_{(l-1)l}|(x_1)^{r-2} \,\,\,\,\,\,\,\,\,\  {\rm by } \,\ (\ref{eq234})  \nonumber \\
&\ge& x_{l-1}^2 \sum_{i_1i_2 \cdots i_{r-2}\in E_{(l-1)l}} x_{i_1}x_{i_2} \cdots x_{i_{r-2}} \nonumber \\
&=& x_{l-1}^2 \lambda(E_{(l-1)l}, \vec{x}).
\end{eqnarray}
This proves part (a).

(b)  Let  ${\vec x}=(x_1, x_2, \ldots, x_l)$ be an optimal weighting for $G$. Since $G$ is left-compressed, by Remark \ref{r1}(a), $x_1\ge x_2 \ge \cdots \ge x_l \ge 0$. If $x_{l-1}=0$, since $G$ does not contain the clique $[l-1]^{(r)}$, then the conclusion holds obviously, so we assume that $x_{l-1}>0$.
 We  add all edges in $[l-1]^{r}-E(G)$ to $G$ and get a new $r$-graph $H$. Observe that  the new $r$-graph $H$ is still left-compressed,  satisfies $|[l-2]^{(r-1)} \backslash E_l| \ge 2^{r-3}|E_{(l-1)l}|$ and contains the clique  $[l-1]^{(r)}$. So by part (a), $\lambda(H)= \lambda([l-1]^{(r)}).$ On the other hand,  $$\lambda(G)=\lambda(G,\vec {x})<\lambda(H,\vec {x})\le \lambda(H).$$
Therefore, $\lambda(G)<\lambda([l-1]^{(r)}).$ This proves part (b).

 \qed

We are now ready to prove Theorem \ref{Corollary 2}.

\bigskip

{\em Proof of Theorem \ref{Corollary 2}.} 
(a) Let $m$ and $l$ be positive integers satisfying ${l-1 \choose r} \le m \le {l-1 \choose r} + {l-2 \choose r-1}-(2^{r-3}-1)({l-2 \choose r-2}-1)$.  Let $G$ be an $r$-graph on $l$ vertices with $m$ edges and a clique of order $l-1$ such that $\lambda(G)=\lambda_{(m,l-1)}^{r}$.
Applying Lemma \ref{lemma1}, we can assume that $G$ is left-compressed and contains the clique $[l-1]^{(r)}$. By Theorem \ref{theorem 2a}, it is sufficient to show that $|[l-2]^{(r-1)} \backslash E_l| \ge 2^{r-3}|E_{(l-1)l}|$.  If not, then  $|[l-2]^{(r-1)} \backslash E_l| < 2^{r-3}|E_{(l-1)l}|$. Since $G$  contains the clique $[l-1]^{(r)}$, then
\begin{eqnarray*}
m&=&{l-1 \choose r}+{l-2 \choose r-1}-|[l-2]^{(r-1)} \backslash E_l|+|E_{(l-1)l}| \\
&>& {l-1 \choose r}+{l-2 \choose r-1}-(2^{r-3}-1)|E_{(l-1)l}| \\
&\ge &{l-1 \choose r}+{l-2 \choose r-1}-(2^{r-3}-1)({l-2 \choose r-2}-1)
\end{eqnarray*}
since $|E_{(l-1)l}|\le {l-2 \choose r-2}-1$. (If $|E_{(l-1)l}|= {l-2 \choose r-2}$, then $E=[l]^{(r)}$ since $G$ is left-compressed and $m={l \choose r}$ which is a contradiction.) This proves part (a) of Theorem \ref{Corollary 2}.

(b) Let $G$ be a $3$-graph with $m$ edges without containing a clique of order $l-1$ such that $\lambda(G)=\lambda_{(m,l-1)}^{3-}$. Then by Lemma \ref{lemmanoclique}, we can assume that $G$ is left-compressed with vertex set $[l]$.
Let ${\vec x}=(x_1, x_2, \ldots, x_l)$ be an optimal weighting of $G$ satisfying $x_1\ge x_2\ge \cdots \ge x_l$. If $x_l$=0, then $\lambda(G)<\lambda([l-1]^{(3)})$ and the conclusion follows. So we assume that $x_l>0$. Now we will use  the following result which is proved in \cite{PZZ}. 

\begin{lemma} \label{lemmaclique} (see \cite{PZZ})
Let $m$ and $l$  be  positive integers satisfying ${l-1 \choose 3} \le m \le {l-1 \choose 3} + {l-2 \choose 2}-(l-2)$. Let $G$ be a left-compressed 3-graph on the vertex set $[l]$ with $m$ edges and without containing a clique of order $l-1$  such that $\lambda(G)=\lambda_{(m,l-1)}^{3-}$. Let ${\vec x}$ be an optimal weighting for $G$ with $l$ positive weights. Then $\lambda(G)< \lambda([l-1]^{(3)}$ or
$$\vert [l-1]^{(3)}\setminus E\vert \le l-2.$$
\end{lemma}

Let $H$ be obtained by adding all triples in $[l-1]^{(3)}\setminus E(G)$ to $G$. By Lemma \ref{lemmaclique}, there are at most $l-2$ such triples. Therefore, $H$ is a $3$-graph with at most ${l-1 \choose 3} + {l-2 \choose 2}$ edges and containing $[l-1]^{(3)}$. By Theorem \ref{theorem 1}, $\lambda(H)=\lambda([l-1]^{(3)})$. Since each $x_i>0$ and $G$ does not contain the clique $[l-1]^{(3)}$, then 
$$\lambda(G)=\lambda(G, {\vec x})<\lambda(H, {\vec x})\le \lambda(H)=\lambda([l-1]^{(3)})$$
which proves part (b) of Theorem \ref{Corollary 2}. \qed

\begin{theo} \label{theorem 3} Let  $l$ be a positive integer.
Let $G$ be a left-compressed $r$-graph on the vertex set $[l]$ such that there is a one-to-one function $f$ from $E_{(l-1)l}$ to $[l-2]^{(r-1)} \backslash E_l$ satisfying the condition:  for $i_1i_2\cdots i_{r-2} \in E_{(l-1)l}$,
  $f(i_1i_2\cdots i_{r-2})=j_1 j_2\cdots j_{r-1}$ satisfies $j_{k}\le i_{k+1}$ for all $1\le k\le r-3$, where $i_1\le i_2\le \cdots \le i_{r-2}$ and $j_1 \le j_2\le \cdots \le j_{r-1}$.

(a) If $G$ contains $[l-1]^{(r)}$, then  $\lambda(G) = \lambda([l-1]^{(r)})$.

(b) If $G$  does not contain $[l-1]^{(r)}$, then  $\lambda(G) <\lambda([l-1]^{(r)})$.
\end{theo}

\bigskip

{\em Proof.} (a)
Let $G$ contain the clique $[l-1]^{(r)}$. 
Let  ${\vec x}=(x_1, x_2, \ldots, x_l)$ be an optimal weighting for $G$. Now we proceed to show that $\lambda(G) \le \lambda([l-1]^{(r)})$.

Consider a new weighting for $G$, ${\vec z}=(z_1, z_2, \ldots, z_l)$ given by $z_i=x_i$ for $i\neq l-1, l$, $z_{l-1}=0$ and $z_l=x_{l-1}+x_l$. By Lemma \ref{LemmaTal5}(a), $\lambda(E_{l-1}, \vec{x})=\lambda(E_{l}, \vec{x})$, so
\begin{eqnarray}\label{eq10a}
\lambda(G,\vec {z})- \lambda(G,\vec {x})
&=&x_{l-1}(\lambda(E_{l}, \vec{x})-x_{l-1}\lambda(E_{(l-1)l}, \vec{x}))\nonumber \\ &&-x_{l-1}(\lambda(E_{l-1}, \vec{x})-x_{l}\lambda(E_{(l-1)l}, \vec{x}))-x_{l-1}x_l\lambda(E_{(l-1)l}, \vec{x}))\nonumber \\
&=&x_{l-1}(\lambda(E_{l}, \vec{x})-\lambda(E_{l-1}, \vec{x}))-x_{l-1}^2\lambda(E_{(l-1)l}, \vec{x})) \nonumber \\
&=& -x_{l-1}^2\lambda(E_{(l-1)l}, \vec{x}).
\end{eqnarray}
We will show that there exists a set of edges $F\subset \{1, ..., l-2, l\}^{(r)}\setminus E$ satisfying
\begin{equation}\label{eq11a}
\lambda(F,\vec {z})\ge x_{l-1}^{2} \lambda(E_{(l-1)l}, \vec{x}).
\end{equation}

Then using (\ref{eq10a}) and (\ref{eq11a}), the $r$-graph $G^{*}=([k], E^{*})$, where $E^{*}=E\cup F$, satisfies $\lambda(G^{*}, \vec {z}))\ge \lambda(G)$. Since $\vec {z}$ has only $l-1$ positive weights, then $\lambda(G^{*}, \vec {z}))\le \lambda([l-1]^{(r)})$, and consequently,
$$\lambda(G)\le \lambda([l-1]^{(r)}).$$
Let $F$ consist of all $j_1 j_2\cdots j_{r-1}l$, where $j_1 j_2\cdots j_{r-1}\in f(E_{(l-1)l})$. 
Then $$\lambda(F,\vec {z}) =(x_{l-1}+x_l)\sum_{j_1 j_2\cdots j_{r-1}\in f(E_{(l-1)l})} x_{j_1} \cdots x_{j_{r-1}}.$$
 Recall that $f$ is a one-to-one function and  for each element $i_1i_2 \cdots i_{r-2}$ of $E_{(l-1)l}$, $f(i_1i_2 \cdots i_{r-2})=j_1 j_2\cdots j_{r-1}$ satisfies $j_{k}\le i_{k+1}$ for all $1\le k\le r-3$.
Combining with Lemma \ref{lemmaadd}, we get
\begin{eqnarray}\label{eq132a}
\lambda(F,\vec {z})
&\ge& x_{l-1}^2 \sum_{i_1i_2 \cdots i_{r-2}\in E_{(l-1)l}} x_1x_{i_2} \cdots x_{i_{r-2}}  \nonumber \\
&\ge& x_{l-1}^2 \sum_{i_1i_2 \cdots i_{r-2}\in E_{(l-1)l}} x_{i_1}x_{i_2} \cdots x_{i_{r-2}}  \nonumber \\
&=& x_{l-1}^2 \lambda(E_{(l-1)l}, \vec{x}).
\end{eqnarray}
This proves (a). 

(b) Let  ${\vec x}=(x_1, x_2, \ldots, x_l)$ be an optimal weighting for $G$. Since $G$ is left-compressed, by Remark \ref{r1}(a), $x_1\ge x_2 \ge \cdots \ge x_l \ge 0$. If $x_{l-1}=0$, since $G$ does not contain the clique $[l-1]^{(r)}$, then the conclusion holds obviously, so we assume that $x_{l-1}>0$.
 We  add all edges in $[l-1]^{r}-E(G)$ to $G$ and get a new $r$-graph $H$. Observe that  the new $r$-graph $H$ is still left-compressed,  contains the clique  $[l-1]^{(r)}$, and still satisfies the condition that there is a one-to-one function $f$ from $E_{(l-1)l}$ to $[l-2]^{(r-1)} \backslash E_l$ such that for $i_1i_2\cdots i_{r-2} \in E_{(l-1)l}$, $f(i_1i_2\cdots i_{r-2})=j_1 j_2\cdots j_{r-1}$ satisfies $j_{k}\le i_{k+1}$ for all $1\le k\le r-3$.
So by part (a), $\lambda(H)= \lambda([l-1]^{(r)}).$ On the other hand,  $$\lambda(G)=\lambda(G,\vec {x})<\lambda(H,\vec {x})\le \lambda(H).$$
Therefore, $\lambda(G)<\lambda([l-1]^{(r)}).$ This proves part (b).
\qed


\begin{remark} \label{remark1}
When $r=3$, Theorem \ref{theorem 3}(a) implies Theorem \ref{theorem 1} and Theorem \ref{theorem 3}(b) implies  Theorem \ref{Corollary 2}(b).
\end{remark}

{\em Proof.}
Let $G$ be a $3$-graph with $m$ edges containing a clique of order $l-1$, where ${l-1 \choose 3} \le m \le {l-1 \choose 3} + {l-2 \choose 2}$. By Lemma \ref{lemmaclique1}, we can assume that $G$  is left-compressed and  on the vertex set $[l]$.   If $\vert E_{(l-1)l}\vert > \vert [l-2]^{(2)}\setminus E_l\vert$, then
$$m\ge {l-1 \choose 3} + {l-2 \choose 2}-\vert [l-2]^{(2)}\setminus E_l\vert+\vert E_{(l-1)l}\vert>{l-1 \choose 3} + {l-2 \choose 2}$$
which is a contradiction. So $\vert E_{(l-1)l}\vert \le \vert [l-2]^{(2)}\setminus E_l\vert$. Therefore, there is a one-to-one function $f$ from $E_{(l-1)l}$ to $[l-2]^{(2)} \backslash E_l$  and $f$ automatically satisfies the condition in Theorem \ref{theorem 3}. Applying Theorem \ref{theorem 3}, we have $\lambda(G) = \lambda([l-1]^{(3)})$. 

Applying Lemma \ref{lemmaclique}, we can similarly show that Theorem \ref{theorem 3}(b) implies  Theorem \ref{Corollary 2}(b).
\qed

In our results below, note that it does not matter how many vertices we are allowed to use.

\begin{theo} \label{theorem3.9}
Let $G$ be an $r$-graph containing a  clique of order $l-1$ with $m$ edges. If $m \le {l-1 \choose r} + 2(l-r)$, then $\lambda(G)=\lambda([l-1]^{(r)})$.
\end{theo}

{\em Proof.}
Let $G$ be an $r$-graph containing a clique of order $l-1$ with $m$ edges such that $\lambda(G)=\lambda_{(m,l-1)}^{(r)}$. Clearly $\lambda(G)\ge\lambda([l-1]^{(r)})$. Next we show that $\lambda(G)\le\lambda([l-1]^{(r)})$.  Since $\lambda_{(m,l-1)}^{(r)}$ does not decrease as $m$ increases, it is sufficient to show  the case that $m={l-1 \choose r} + 2(l-r)$.   Based on Lemma \ref{lemma1}, we may assume that $G$ is  left compressed and the optimal weighting ${\vec x}=(x_1, x_2, \ldots, x_n)$ of $G$ satisfying $x_i\ge x_j$ when $i<j$. Note that $x_{l+1}=0$. Otherwise, then  by Lemma \ref{LemmaTal5},  $G$ contains edge $12 \ldots (r-2)l(l+1)$. Since $G$ is left-compressed,  $G$  contains all edges $12 \ldots (r-2)i(l+1)$  for all $i$, where $r-1 \le i \le l$ and $12 \ldots (r-2)jl$ for all $j$, where $r-1 \le j \le l-1$. Then $m\ge {l-1 \choose r} + 2(l-r)+3$ which is a contradiction. If $x_l=0$, then $\lambda(G)\le\lambda([l-1]^{(r)})$ and we are done. If $x_l>0$, then  by Lemma \ref{LemmaTal5},  $G$ contains edge $12 \ldots (r-2)(l-1)l$, it should contain all edges $12 \ldots (r-2)il$ for all $i$, where $r-1 \le i \le l-1$. Note that $E_{(l-1)l}=\{12\ldots(r-2)\}.$ Otherwise,  $12\cdots (r-3)(r-1)(l-1)l \in E$ and $G$ contains all edges $12 \cdots (r-3)(r-1)il$ for all $i$, where $r \le i \le l-1$. So $m\ge {l-1 \choose r} + 2(l-r)+1$ which is a contradiction. So
 we can assume that $G$ is on the vertex set $[l]$ and $E_{(l-1)l}=\{12 \cdots (r-2)\}$.

Consider a new weighting for $G$, ${\vec z}=(z_1, z_2, \ldots, z_l)$ given by $z_i=x_i$ for $i\neq l-1, l$, $z_{l-1}=0$ and $z_l=x_{l-1}+x_l$. By Lemma \ref{LemmaTal5}(a), $\lambda(E_{l-1}, \vec{x})=\lambda(E_{l}, \vec{x})$, so
\begin{eqnarray}\label{eq20}
\lambda(G,\vec {z})- \lambda(G,\vec {x})
&=&x_{l-1}(\lambda(E_{l}, \vec{x})-\lambda(E_{l-1}, \vec{x}))-x_{l-1}^2 \lambda(E_{(l-1)l},\vec{x}) \nonumber \\
&=& -x_{l-1}^2 \lambda(E_{(l-1)l},\vec{x}) \nonumber \\
&=&-x_{l-1}^2x_1x_2 \cdots x_{r-2}.
\end{eqnarray}

Let $F$ consist of those edges in $\{1, ..., l-2, l\}^{(r)}\setminus E$ containing the vertex $l$. Then clearly
\begin{equation}\label{eq232}
\lambda(F,\vec {z})\ge (x_{l-1}+x_l)x_2 \cdots x_{r-2}x_{l-1}^2
\end{equation}
since  $23\cdots (r-2)(r-1)(l-2)l$  is in $F$. Otherwise, $23\cdots (r-2)(r-1)(l-2)l\in E$. Since $G$ is left-compressed, then $23\cdots (r-2)(r-1)il\in E$  for all $r\le i\le l-2$, $12\cdots (r-2)jl\in E$ for all $r-1\le j\le l-2$ and $13\cdots (r-1)kl\in E$ for all   $r\le k\le l-2$. Recall that $12 \ldots (r-2)il\in E$ for all $i$, where $r-1 \le i \le l-1$.
Then $m\ge {l-1 \choose r}+3(l-r)-1\ge {l-1 \choose r}+2(l-r)+1$ which is a contradiction.


By Lemma \ref{lemmaadd}, we have $x_1 \le x_{l-1}+x_l$. Applying this to (\ref{eq232}), we get 
\begin{equation}\label{eq21}
\lambda(F,\vec {z})\ge x_{l-1}^2x_1x_2 \cdots x_{r-2}.
\end{equation}

Then using (\ref{eq20}) and (\ref{eq21}), the $r$-graph $G^{*}=([l], E^{*})$, where $E^{*}=E\cup F$, satisfies $\lambda(G^{*}, \vec {z}))\ge \lambda(G)$. Since $\vec {z}$ has only $l-1$ positive weights, then $\lambda(G^{*}, \vec {z}))\le \lambda([l-1]^{(r)})$, and consequently,
$$\lambda(G)\le \lambda([l-1]^{(r)}).$$
This completes the proof. \qed


\section{Concluding Remarks}

As we have seen in  Section \ref{sectionconjecture1}, Lemma \ref{lemmaclique1}, Lemma \ref{lemmanoclique}, and Theorem \ref{theorem 2a} for the case $r=3$ refine Theorem \ref{Tal}.   If one can have some results similar to Lemma \ref{lemmaclique1} and Lemma \ref{lemmanoclique} for general $r$, then one can get results similar to Theorem \ref{Tal} for general $r$.

We also remark that in some applications, estimating $\lambda(l,r,m)$ is sufficient. So Theorems \ref{Corollary 2} and \ref{Talr} might still be applicable in some situations though Conjectures \ref{conjecture1} and \ref{conjecture2} cannot be verified in general at this moment.

\bigskip

{\bf Acknowledgments.} We thank an anonymous referee and the editor for  helpful and insightful comments.


\begin{thebibliography}{1}


\bibitem{B2} M.  Budinich, Exact bounds on the order of the maximum clique of a graph, Discrete Appl. Math. 127, 535-543 (2003).

\bibitem{B3} S.  Busygin, A new trust region technique for the maximum weight clique problem,  Discrete Appl. Math. 304(4), 2080-2096 (2006).

\bibitem{FF} P. Frankl and Z. F\"uredi, Extremal problems whose solutions are the blow-ups of the small Witt-designs,  Journal of Combinatorial Theory (A) 52 (1989), 129-147.

\bibitem{FR84} P. Frankl and V. R\"{o}dl,  Hypergraphs do not jump,  Combinatorica 4 (1984), 149-159.

\bibitem{G9}  L. E. Gibbons, D. W. Hearn, P. M. Pardalos, and M. V. Ramana, Continuous characterizations of the maximum clique problem, Math. Oper. Res., 22 (1997), 754-768.

\bibitem{HPZ} G. He, Y. Peng, and C. Zhao, On finding Lagrangians of 3-uniform hypergraphs, Ars Combinatoria (accepted).

\bibitem{MS} T.S. Motzkin and E.G. Straus, Maxima for graphs and a new proof of a theorem of Tur\'an,  Canad. J. Math 17 (1965), 533-540.

\bibitem{mubayi06} D. Mubayi, A hypergraph extension of Tur´an's theorem, J. Combin. Theory
Ser. B 96 (2006), 122-134.

\bibitem{PP} M. Pavan and M. Pelillo, Generalizing the motzkin-straus theorem to edge-weighted graphs, with applications to image segmentation,  Lecture Notes in Computer Science 2683(2003), 485-500.

\bibitem{PP15} P.M. Pardalos and A.T. Phillips,  A global optimization approach for solving the maximum clique problem,  Int. J. Comput. Math. 33 (1990), 209-216.

\bibitem{PZ} Y. Peng and C. Zhao, A Motzkin-Straus type result for 3-uniform hypergraphs, Graphs and Combinatorics, 2012, DOI: 10.1007/s00373-012-1135-5.

\bibitem{PZ2} Y. Peng and C. Zhao, On Lagrangians of hypergraphs and cliques, Recent Advances in Computer Science and Information Engineering 125 (2012), 7-12.

\bibitem{PZZ} Y. Peng, H. G. Zhu, and C. Zhao, On cliques and Lagrangians of 3-uniform hypergraphs, submitted. http://arxiv.org/abs/1211.6508

\bibitem{BP1} S. Rota Bul\'o, and M. Pelillo,  A continuous characterization of maximal cliques in k-uniform hypergraphs, In Learning and Intellig. Optim.,  Vol.5313 (2008), 220-233.

\bibitem{BP2} S. Rota Bul\'o and M. Pelillo, A generalization of the Motzkin-Straus theorem to hypergraphs, Optim. Letters 3, 2 (2009), 287-295.

\bibitem{RTP20} Rota Bul\'o, S., Torsello, A., and M. Pelillo, A continuous-based approach for partial clique enumeration, Graph-Based Representations Patt. Recogn. 4538 (2007), 61-70.

\bibitem{sidorenko89} A. F. Sidorenko, Solution of a problem of Bollob\'as on 4-graphs,  Mat. Zametki 41 (1987), 433-455.

\bibitem{SS} V. S\'os, and E. G. Straus,  Extremal of functions on graphs with applications to graphs and hypergraphs, J. Combin. Theory Series B 63 (1982), 189-207.


\bibitem{T} J. Talbot, Lagrangians of hypergraphs,  Combinatorics, Probability \& Computing 11 (2002), 199-216.

\bibitem{TPZZ1} Q. S. Tang, Y. J. Peng, X. D. Zhang, and C. Zhao, Some results on Lagrangians of hypergraphs, submitted. http://arxiv.org/abs/1211.7057

\bibitem{TPZZ2} Q. S. Tang, Y. J. Peng, X. D. Zhang, and C. Zhao, On Frankl and F\"{u}redi's conjecture for $3$-uniform hypergraphs, submitted. http://arxiv.org/abs/1211.7056

\bibitem{Turan}P. Tur\'an, On an extremal problem in graph theory(in Hungarian),  Mat. Fiz. Lapok 48 (1941), 436-452.

\bibitem{wilf} H.S. Wilf, Spectral bounds for the clique and independence number of graphs,  J. Combin. Theory Ser. B 40 (1986), 113-117.


\end{thebibliography}
\end{document}